\documentclass[draft]{amsart}

\usepackage{amssymb,amscd}
\usepackage{verbatim}
\usepackage[T1]{fontenc}

\numberwithin{equation}{section}
\newtheorem{thm}[equation]{Theorem} 
\newtheorem{lem}[equation]{Lemma}
\newtheorem{prop}[equation]{Proposition} 
\newtheorem{cor}[equation]{Corollary}

\newtheorem{exthm}[equation]{Example}

\theoremstyle{definition} 

\newtheorem{rem}[equation]{Remark} 

\def \fm{\mathfrak m}

\def\aa{a}
\def\bb{b}

\def\II{I\kern-3pt I}

\def\preisomto{\vbox{\hbox to 13pt{\hfill$\hskip -0.5pt
                 \sim$\hfill}\nointerlineskip\vskip -0.5pt \hbox to
                 13pt{\rightarrowfill}}}

\def\prelongisomto{\vbox{\hbox to
                16pt{\hfill$\sim$\hfill}\nointerlineskip\vskip -0.5pt
                \hbox to 16pt{\rightarrowfill}}}

\newcommand{\doublelongrightarrow}{\longrightarrow \kern-14pt
\longrightarrow}

\def\trait{\hbox to 2mm{\hrulefill}}

\def\G-rank{\operatorname{Grank}}

\def\Hom{\operatorname {Hom}}
\def\End{\operatorname {End}}

\def\gr{\operatorname {gr}}

\def\Ext{\operatorname {Ext}}

\def \Rees{\operatorname{Rees}}

\def\gr{\operatorname{gr} }

\begin{document} 

\title{Algebras without noetherian filtrations}

\author{J. T. Stafford} 
\address{Department of Mathematics, University of Michigan,
Ann Arbor, MI 48109-1109, USA.}
\email{jts@math.lsa.umich.edu}

\author{J. J. Zhang}
\address{Department of Mathematics, Box 354350, University of
Washington, Seattle, WA 98195, USA.}
\email{zhang@math.washington.edu}

\thanks{Both authors were supported in part by the NSF. The second author
was  also supported by the Royalty Research Fund of the University
of Washington}

\begin{abstract}
We provide examples of finitely generated noetherian  PI algebras for which
there is no finite dimensional filtration with a noetherian associated graded
ring; thus we answer 
negatively a question of Lorenz \cite[p.436]{Lo}.
\end{abstract}
\maketitle


\section{Introduction}\label{intro}

In this paper all algebras will be defined over a fixed base field $k$.
Let $\Gamma$ be an  $\mathbb N$-filtration of a
$k$-algebra $A$; thus, $\Gamma = \{\Gamma_i : i\in \mathbb N\}$
is an ascending chain of $k$-subspaces of $A$ satisfying $1\in \Gamma_0$, 
$ \bigcup_{j\in \mathbb N} \Gamma_j = A$ 
and $\Gamma_i\Gamma_j\subseteq \Gamma_{i+j}$ for all
$i,j$. The filtration is  {\it finite } if $\dim_k \Gamma_i<\infty$ for all
$i\in \mathbb N$ and {\it standard} if $\Gamma_0=k$ and $\Gamma_i=\Gamma_1^i$
for all $i\geq 2$.  We say that $\Gamma$ is a {\it  
(left) noetherian filtration} if the associated graded 
ring $\gr A= \gr_{\Gamma}A = \bigoplus_{i} \Gamma_i/\Gamma_{i-1}$
 is (left) noetherian. The algebra $A$ is  {\it affine} if it is finitely
generated as a $k$-algebra.

If an algebra $A$ has a left noetherian $\mathbb N$-filtration,
 then a standard technique is to pull results back from $\gr A$ to
$A$ since theorems are typically easier to prove in the graded ring
$\gr A$ than in $A$. 
For related reasons Lorenz asked in \cite[p.436]{Lo}
and \cite[Question III.4.2]{Lo2}
whether every left noetherian affine algebra $R$, satisfying 
a polynomial identity (PI), admits a left noetherian, standard finite 
$\mathbb N$-filtration. This question has 
been raised again (without the ``standard'' hypothesis)
 in  \cite[Question 6.16]{YZ} because of its importance 
for dualizing complexes and homological questions:
 if Lorenz's question were to have a positive answer then 
$R$ would have a dualizing complex \cite[Corollary 6.9]{YZ}
and  every noetherian affine PI Hopf algebra 
would have finite injective dimension \cite{WZ}.

The aim of this note is to answer these questions by 
providing a class of noetherian affine 
PI algebras which do not admit any noetherian finite $\mathbb N$-filtration.
The basic technique is provided by the following theorem (see
Section~\ref{sec3}).

\begin{thm}\label{thm1.1} Suppose that $I$ and $J_1\subset J_2$
are ideals of an algebra $R$  such that $J_2/J_1$ is free of rank $s$ 
as a left $R/I$-module and free of rank $t$ as a right $R/I$-module.
If $s<t$, then there is no finite 
$\mathbb N$-filtration $\Gamma$ of $R$ such that
$\gr_\Gamma R$ is  left noetherian.
\end{thm}

A variant of this theorem also holds if one replaces ``free of rank $x$'' by
``of Goldie rank $x$.'' See Theorem~\ref{thm1.11} for the
details.

A  simple example satisfying the hypotheses (and conclusion) of the theorem
is given by the ring
\begin{equation}\label{ex2.1}
R=\left\{ \pmatrix f(x)& g(x) \\ 0 & f(x^2)\endpmatrix: f,g\in
k[x]\right\}\;\subset\; M_2(k[x]).
\end{equation}
Here, one takes $J_2=I$ to be the ideal of strictly upper triangular matrices 
and $J_1=0$. See Section~\ref{sec4} for more details. 

The most important case of Theorem~\ref{thm1.1}
is when $\Gamma$ is a standard. 
However, we emphasize that the theorem holds 
for any filtration satisfying the earlier definitions.

The analogue of Theorem~\ref{thm1.1} also holds for $\fm$-adic filtrations 
(where $\fm$ is the Jacobson radical of a local algebra)
 and for that reason we prove the result for rings with a
Zariskian filtration (see Theorem~\ref{thm3.3}).
In particular, we  provide an example of a local prime noetherian PI 
ring for which the
Jacobson radical does not satisfy the strong AR property. This is given in
Section~\ref{sec4} where the reader may find further applications of the main
theorem.


\bigskip
\section{Proof of Theorem~\ref{thm1.1}}\label{sec3}

\bigskip

Since we are also interested in $\fm$-adic filtrations
of local rings as well as ascending filtrations,
 we will prove our main result for  $\mathbb Z$-filtrations.
First we review some basic facts about filtrations  
from \cite[Chapter 6]{KL} and \cite{LV}. 

Suppose that $A$ is a $k$-algebra. For the purposes of this paper
a  {\it  filtration} (or 
more strictly, an exhaustive separated finite $\mathbb Z$-filtration) 
of $A$ is an ascending chain   of subspaces
$\Gamma=\{\Gamma_i\subseteq \Gamma_{i+1}\;|\; i\in {\mathbb Z}\}$ of $A$,
 satisfying:

\vskip 2pt
\begin{enumerate}
\item[(1)] $1\in \Gamma_0$  and  $\Gamma_i \Gamma_j\subseteq 
\Gamma_{i+j}$ 
for all $i,j\in {\mathbb Z}$;
\item[(2)]
 $\Gamma$ is {\it finite} in the sense that 
$\dim \Gamma_i/\Gamma_{i-1}<\infty$  
for all $i\in {\mathbb Z}$;
\item[(3)] $\Gamma$ is {\it exhaustive} in the sense that 
 $A=\bigcup_{i\in {\mathbb Z}} \Gamma_i$ and {\it separated} 
in the sense that $\bigcap_{i\in {\mathbb Z}} 
\Gamma_i=0$.
\end{enumerate}

\vskip 2pt
\noindent
The {\it Rees ring} associated to $\Gamma$ is defined to be
$\Rees A=\Rees_\Gamma A=\bigoplus_{i\in {\mathbb Z}} \Gamma_i$
and the {\it associated graded ring} is 
$\gr A=\gr_{\Gamma} A=\bigoplus_{i\in {\mathbb Z}} \Gamma_i/\Gamma_{i-1}.$
Write
$J(A)$ for the  Jacobson radical of  $A$. 
Following \cite{LV}, 
a filtered algebra $A$ is called {\it (left) Zariskian}  if
 
\vskip 2pt
(Zar1) $\Gamma_{-1}\subseteq J(\Gamma_0)$;

(Zar2) $\Rees_F A$ is (left) noetherian.

\vskip 2pt
Fix a filtration $\Gamma$ of the algebra $A$ and a left $A$-module $M$. 
The concept of a $\mathbb Z$-filtration (again, finite, exhaustive
and separated)
 $\Lambda =\{\Lambda_i : i\in \mathbb Z\}$
of  $M$ is defined analogously;
one simply replaces (1) in the above definition by

\vskip 2pt
(1$'$) $\Gamma_i\Lambda_j \subseteq \Lambda_{i+j}$ 
for all $i,j\in {\mathbb Z}$.

\vskip 2pt
\noindent
Corresponding to this filtration one has the Rees module
 $\Rees_\Lambda M=\bigoplus \Lambda_i$ and 
associated graded module $\gr_\Lambda M=\bigoplus \Lambda_i/\Lambda_{i-1}$.
We say that $\Lambda$ is a {\it good filtration} if there exist $\{m_i\in 
\Lambda_{d_i} : 1\leq i\leq r<\infty\}$ such that $\Lambda_n=\sum_{i=1}^r
\Gamma_{n-d_i}m_i$ for all $n\in\mathbb Z$.
Two filtrations $\Lambda$ and $\Lambda'$ of $M$ are 
{\it equivalent}, written $\Lambda\sim \Lambda'$, if
there is an integer $q$ such that
$\Lambda_i\subset \Lambda'_{i+q}$ and $
\Lambda'_i\subset \Lambda_{i+q} $ for all $
i\in {\mathbb Z} 
.$


The {\it Hilbert
function} of $M$ with respect to $\Lambda$ is defined
to be 
$$H_{M,\Lambda}(n)= \dim \Lambda_n/\Lambda_{-n}, 
\quad {\text{for all}}\
n\in {\mathbb Z}.$$
Let $H$ and $H'$ be two (Hilbert) functions $\mathbb N\to \mathbb N$. 
We say that the growth of $H$ is at most the growth of $H'$, and write
$H\leq H'$, if there is
an integer $q$ such that $H(n)\leq H'(n+q)$ for all 
$n\gg 0$.
We say $H$ and $H'$ are 
{\it equivalent}, written $H\sim H'$,
 if both $H\leq H'$ and $H'\leq H$ hold.
  
Since we require filtrations to be separated, they need not induce filtrations
on factor modules. However, for Zariskian filtrations this is not a problem:

\begin{lem}\label{lem3.1} Suppose that $\Gamma=\{\Gamma_i\}$ is a 
filtration of an algebra $A$ and let $M$ be a left $A$-module
with  a good filtration $\Lambda$.
\begin{enumerate}
\item[(1)]
If $\Lambda'$ is another filtration of $M$,   there exists an 
integer $q$ such that $\Lambda_i\subset \Lambda'_{i+q}$ for all $i$.
 Thus, if $\Lambda'$ is good, then $\Lambda$ and $\Lambda'$ are
equivalent.
\item[(2)]
If $\Lambda'$ and $\Lambda''$ are equivalent filtrations of $M$,
then $H_{M,\Lambda'}$ and $H_{M,\Lambda''}$ are 
equivalent. 
\item[(3)]
Assume that $\Gamma$ is left Zariskian.
 If $N$ is a submodule of $M$ then $\Lambda$ 
induces good filtrations on $N$ and $M/N$. In particular,
 $\Gamma$ induces a left Zariskian filtration on
every factor ring of $A$. 
\end{enumerate}
\end{lem}

\begin{proof} (1) and (2) follow from the definitions while (3) follows 
from \cite[Theorem~3.3]{LV} or \cite[Theorem II.2.1.2]{LV2}. 
\end{proof}

The key observation in this paper is given by the next proposition. 
Since our  main theorem will come in three slightly different forms, 
this result will also have three slightly different cases.

\begin{prop}\label{prop3.2} Let $A$ be an algebra with a filtration
$\Gamma$ such that the growth of $H_{A,\Gamma}$ is subexponential. 
Let $M$ be an $A$-bimodule such that the left module
$_AM$ is free of rank $s$ and the right module $M_A$ is free of rank $t$. 
Suppose that there exists a filtration $\Lambda$ on $M$ such that 
$\Lambda$ is a good filtration  of ${}_AM$
and a filtration of $M_A$. 
\begin{enumerate}
\item[(1)] If $\Lambda $ is also a good filtration of $M_A$ then $t=s$.
\item[(2)] If $\Gamma_n=0$ for $n\ll 0$ then $s\geq t$
\item[(3)] If $\Gamma_n=A$ for $n\gg 0$ then $s\leq t$.
\end{enumerate}
\end{prop}

\begin{proof}  The beginning of the proof is the same in all three cases.
Note that ${}_AA$ has a good filtration, simply because $\Gamma_i=\Gamma_i1$ for
all $i$. Thus ${}_AA^{(s)}$ also has a good filtration and 
so this induces a 
  good filtration $\Lambda'$ on $M$ such that $\Lambda'_i\cong \Gamma_i^{(s)}$
  for all $i$.
   By Lemma~\ref{lem3.1}(1), there exists $q_1$ such that 
 $\Lambda_i\subseteq \Lambda'_{i+q_1}$ for all $i$.

Similarly,   the right $A$-module structure provides
an induced good filtration $\Lambda''$ of $M_A$ such that 
$\Lambda''_i\cong \Gamma_i^{(t)}$  for all $i$.
 By Lemma~\ref{lem3.1}(1), there exists $q_2$ such that 
 $\Lambda''_i\subseteq \Lambda_{i+q_2}$ for all $i$.
Hence, for $p=q_1+q_2$, 
\begin{equation}\label{formula}
\Lambda''_i\subseteq \Lambda'_{i+p} \quad {\rm for\ all \ } i\in \mathbb Z.
\end{equation}

We now consider the three cases separately. Under assumption (2),
we see that 
$\Lambda''_i/\Lambda''_{-i}=\Lambda''_i \hookrightarrow \Lambda'_{i+p} =
\Lambda'_{i+p}/\Lambda'_{-(i+p)}$ for $i\gg 0$.
Thus, $H_{M,\Lambda''}(n)\leq H_{M,\Lambda'}(n+p)$.
But, by construction, $H_{M,\Lambda'}=sH_{A,\Gamma}$
and $H_{M,\Lambda''} = tH_{A,\Gamma}$ and so 
$$  H_{A,\Gamma}(n)\leq \frac{s}{t} H_{A,\Gamma}(n+p)
\quad {\rm  for \ all \ } n\in \mathbb N.$$
Since $H_{A,\Gamma}(n)$ grows subexponentially this forces 
$t\leq s$.

Under assumption (3) we have 
$\Lambda''_{i+p}/\Lambda''_{-(i+p)}=A/\Lambda''_{-i-p} \twoheadrightarrow
A/\Lambda'_{-i} =\Lambda'_i/\Lambda'_{-i}$ for $i\gg 0$. 
Thus,  $H_{M,\Lambda'}(i)\leq H_{M,\Lambda''}(i+p)$ and
repeating the analysis of the last paragraph shows that $s\leq t$.

Finally, assume that (1) holds. In this case, 
$\Lambda''\sim\Lambda$ and so $\Lambda''\sim\Lambda'$.
Thus there exists $p$ such that $\Lambda'_i\subseteq \Lambda_{i+p}''
\subseteq \Lambda'_{i+2p}$.
Now a minor variant of the penultimate paragraph shows that
$t\leq s$ and hence, by symmetry, that $s=t$.
\end{proof}

\begin{thm}\label{thm3.3} Suppose that $R$ is an algebra with ideals
$I$ and $J_1\subset J_2$ such that $J_2/J_1$ is free of rank $s$ 
as a left $R/I$-module and free of rank $t$ as a right $R/I$-module.
If $s<t$, then there is no filtration of $R$ such that
$R$ is both left and right Zariskian.
\end{thm}

\begin{proof} Let $\widetilde{\Gamma}$ be a 
left and right Zariskian filtration of $R$. 
Then Lemma~\ref{lem3.1}(3) implies that 
the induced filtration $\Gamma$ on $A=R/I$ is left and right Zariskian 
and so  $\gr_{\Gamma} A$ is 
left (and right) noetherian. By  \cite[Remark after 1.2]{SZ} the growth
of $H_{A,\Gamma}$ is subexponential.
 By Lemma~\ref{lem3.1}(3), 
  $\Gamma$ induces a good filtration on $J_2/J_1$
as a left and a right $A$-module, which contradicts to
Proposition~\ref{prop3.2}(1).
\end{proof}

A curious feature of this result is that $R$ could be left or right Zariskian;
it just cannot be both (see Corollary~\ref{cor4.fun}). 
However, with a little more information one can determine which side goes
wrong.

\begin{thm}\label{thm3.4} Suppose that $R$ is an algebra with ideals
$I$ and $J_1\subset J_2$ such that $J_2/J_1$ is free of rank $s$ 
as a left $R/I$-module and free of rank $t$ as a right $R/I$-module.
Assume that  $s<t$ and that $\Gamma$ is a filtration of $R$.
\begin{enumerate}
\item[(1)] If $\Gamma_i=0$ for $i\ll 0$, then $\Gamma$ is not left Zariskian.
\item[(2)] If $\Gamma_i=A$ for $i\gg 0$, then $\Gamma$ is not right Zariskian.
\end{enumerate}
\end{thm}

\begin{proof}  For (1),  repeat the proof of Theorem~\ref{thm3.3}, but with 
  Proposition~\ref{prop3.2}(1) 
 replaced by  Proposition~\ref{prop3.2}(2). For (2), use the proof of 
 Theorem~\ref{thm3.3}, applied to the opposite ring $R^{\rm op}$, with 
  Proposition~\ref{prop3.2}(1) 
 replaced by  Proposition~\ref{prop3.2}(3).
\end{proof}

We are now ready to prove Theorem~\ref{thm1.1} from the introduction.

\vskip 3pt 
\noindent
{\it Proof of Theorem~\ref{thm1.1}.}
It is easy to see that an $\mathbb N$-filtration $\Gamma$ of a ring $R$ is 
left Zariskian if and only if $\gr_\Gamma$ is left noetherian (see, for example,
\cite[Proposition~II.1.2.3]{LV2}). Thus, 
 Theorem~\ref{thm1.1} is a special case of
Theorem~\ref{thm3.4}(1).
\qed
\vskip 3pt

 The analogue of Theorem~\ref{thm1.1} for complete
local rings also holds with much the
same proof. To state the result, we need some definitions.
Let  $R$ be a semilocal  
algebra with Jacobson radical $J(R)=\fm$ and 
assume that $\dim_kR/\fm <\infty$. A filtration $\Gamma$ of 
 $R$ is called a {\it weak  adic filtration} if 
it satisfies $\Gamma_n=A$ for all $n\geq 0$ and 
 $\Gamma_{-1}\subseteq \fm$. (We still require  the 
filtration to be finite, separated and exhaustive.)

\begin{cor}\label{cor3.7} Let $R$ be a complete semilocal  right noetherian 
algebra
with Jacobson radical $\fm$. Suppose that $R$ has ideals
$I$ and $J_1\subset J_2$ such that $J_2/J_1$ is free of rank $s$ 
as a left $R/I$-module and free of rank $t$ as a right $R/I$-module.
If $s<t$, then there is no 
weak adic filtration $\Gamma$ such that 
$\gr_{\Gamma} R$ is right noetherian. 
\end{cor}

\begin{proof}  
Suppose that such a filtration $\Gamma $ exists. 
By  \cite[Corollary 1.2]{Jat} the filtration $\Gamma$ is  also complete
(in the natural sense that 
Cauchy sequences modulo the $\Gamma_{-i}$ should converge---see
 \cite[Definition I.3.3.2]{LV2}). 
Thus \cite[Proposition~II.2.2.1]{LV2}
implies  that the filtration is right Zariskian and 
the result   follows Theorem~\ref{thm3.4}(2).
\end{proof}

 An analogue of this corollary also holds for non-complete rings, 
although the result is less pleasant since one cannot now assume that 
separated filtrations induce separated filtrations on factor modules.
The result becomes the following: Suppose that $R$ satisfies the 
hypotheses of Corollary~\ref{cor3.7}, but that $R$ is not complete.
Let $\Gamma$ be a weak adic
filtration of $R$ that induces both a filtration 
 on $R/I$ and a good right filtration of $J_2/J_1$. 
If $s<t$, then 
$\gr_{\Gamma} R$ is not  right noetherian. 

\section{A partial generalization}\label{newthm}

Although the results of the last section are sufficient for 
our examples, 
one can give a version of  Theorem~\ref{thm1.1} that works   
without the assumption that $J_2/J_1$ be free,
but at the expense of of assuming that $R/I$ be a prime Goldie ring.
We prove this in this section.
Let $A$ be a prime left Goldie ring with simple artinian ring of fractions
$Q(A)$. If 
 $M$ is a  left $A$-module, then the 
{\it Goldie rank} of $M$ is defined to be the length of 
$Q(A)\otimes_AM$ and written $\G-rank(M)$.

\begin{lem} 
\label{add2}
Suppose that $A$ is a prime Goldie ring with  a left
noetherian 
$\mathbb N$-filtration $\Gamma$. Let $M$ be an  
  $A$-bimodule with a filtration 
    $\Lambda$ that is a good filtration of ${}_AM$ and
  a filtration of $M_A$.  
 If  ${}_AM$ is torsion, then so is $M_A$.
\end{lem}

\begin{proof} Assume that $M_A$ is not torsion. Replacing $M$ 
by  $M^{(n)}$, for some $n$, we may assume that $A_A
\subseteq  M_A$.  
Thus,  Lemma~\ref{lem3.1}(1) implies that 
there exists $q\geq 0$ such that 
$$\dim \Gamma_n \leq \dim \Lambda_{n+q}\cap A\leq \dim \Lambda_{n+q}$$
for all $n$. Since $\Gamma_n=\Lambda_n=0$ for $n\ll 0$, this implies that
$H_{M,\Lambda}\geq H_{A,\Gamma}$.

Now consider ${}_AM$. Since $\Lambda$ is a good filtration, ${}_AM$ is finitely
generated.
We claim, for all $p>0$, that 
$H_{M,\Lambda}\leq \frac{1}{p} H_{A,\Gamma}$.
Once this has been proved, then the last paragraph implies, for some
$x>0$,  that
$H_{A,\Gamma}(n) \leq \frac{1}{p}H_{A,\Gamma}(n+x)$.
This contradicts the fact that $H_{A,\Gamma}$ grows subexponentially
\cite[Remark after 1.2]{SZ} and proves the lemma. 
In  order  to prove the claim we ignore the right-hand structure
of $M$ and so, by induction, it suffices to prove it for a cyclic module
$M=A/I$. Since $_AM$ is torsion and $A$ is Goldie, $I$ contains
a regular element, $a$ say, of $A$. 
We will still write $\Gamma$ for the good filtration 
on any subfactor of ${}_AA$ induced from $\Gamma$.
Now, for any $n$, $M$ is a homomorphic image of 
$L(i)=Aa^{i-1}/Aa^{i}\cong A/Aa$ and so 
$H_{M,\Lambda}\leq H_{L(i), \Gamma}$.
Since Hilbert series are additive on short exact sequences, 
this implies that
$pH_{M,\Lambda}\leq H_{A/Aa^p, \Gamma}
\leq H_{A,\Gamma}$, for any $p\in \mathbb N$. 
\end{proof}


\begin{thm}\label{thm1.11} Suppose that $I$ and $J_1\subset J_2$
are ideals of an algebra $R$  such that  $A=R/I$ is a prime 
Goldie ring. Assume that   $J_2/J_1$ has
Goldie rank $s$ as a left $A$-module and 
Goldie rank $t$ as a right $A$-module, for some 
$s<t$. Then there is no
$\mathbb N$-filtration $\Gamma'$ of $R$ such that
$\gr_{\Gamma'} R$ is  left noetherian.
\end{thm}

\begin{proof}  
Suppose that such a filtration $\Gamma'$  exists.  Let 
$\Gamma$ be the induced filtration on $A$ and $\Lambda$  the induced
filtration on $M=J_2/J_1$. Then $\Lambda$ is a good left filtration and so
 ${}_AM$ is  finitely generated.
The torsion submodule $T$ of $M_A$ is an $A$-bimodule and so we may pass to
$M/T$ without affecting the hypotheses (although $s$ may decrease).
 By Lemma~\ref{add2}
${}_AM$ is also torsion-free.  If $\G-rank(A)=u$, replace $M$ by 
 $M^{( u)}$; thus $A_A^{( t)}\cong X \subseteq M_A $ and  
${}_AM\subseteq Y\cong {_AA}^{( s)}$.
Let $\Theta$, respectively $\Phi$, be the filtrations of $X$ and $Y$ induced
from $\Gamma$.
Since $\Theta$  equals the direct sum $\Gamma^{(t)}$
of $t$ copies of $\Gamma$, certainly $H_{X,\Theta}=tH_{A,\Gamma}$.
Similarly, $H_{Y,\Phi}=sH_{A,\Gamma}$.

Since $\Lambda$ is a good left filtration and $\Theta$ is a good right
filtration, by Lemma~\ref{lem3.1}(1) there exists $q\geq 0$ such that 
$$\dim \Theta_n \leq \dim (X\cap \Lambda_{n+q}) \leq \dim \Lambda_{n+q}$$
and
$$
\dim \Lambda_n \leq \dim ( M\cap \Phi_{n+q}) \leq \dim \Phi_{n+q},$$
for all $n\gg 0$.
Hence, $$tH_{A,\Gamma}(n)=H_{X,\Theta}(n) \leq 
H_{Y,\Phi}(n+2q) = sH_{A,\Gamma}(n+2q),$$
for all $n\gg0$.  Since $H_{A,\Gamma}$ grows subexponentially, this 
forces $t\leq s$.
\end{proof}

With a rather more complicated argument one can prove an analogous 
version of Theorem~\ref{thm3.3} using Goldie ranks.
 However, we do not know how to prove the
analogous version of Theorem~\ref{thm3.4}(2) or Corollary~\ref{cor3.7} 
without extra hypotheses. 

\section{Examples}\label{sec4}

In this section we use the results of the Section~\ref{sec3} to provide
examples of rings without noetherian associated graded rings. 

\begin{exthm}\label{ex4.1} Let 
$$S=\left\{ \pmatrix f(x)& g(x) \\ 0 & f(x^2)\endpmatrix: f, g\in
k[x]\right\}+yM_2(k[x,y])\;\subset\; M_2(k[x,y]).$$ 
Then $S$ is an affine  noetherian prime PI algebra without any 
left noetherian  finite $\mathbb N$-filtrations.
\end{exthm}

\begin{proof} The diagonal matrices of the form
$$\pmatrix
f(x,y) & 0\\ 0 & f(x^2,y) \endpmatrix : f(x,y)\in k[x,y]$$
clearly form a subring $C$  of $S$ isomorphic to $k[x,y]$.
Since $S$ is finitely generated as a left or right $C$-module, it follows that 
 $S$ is a finitely generated noetherian PI algebra.
  It is prime since it contains a nonzero
ideal of the prime ring $M_2(k[x,y])$.

Suppose that  $A\to B$ is a surjective 
homomorphism of algebras and that  $B$ does not admit a left noetherian 
$\mathbb N$-filtration. Then an  immediate consequence of Lemma~\ref{lem3.1}(3)
is that   $A$ also has no such filtration.  Thus, it suffices to prove the 
final assertion for a factor ring, and we use 
\begin{equation}\label{ex2.11}
R=S/yM_2(k[x,y]) \;\cong\; \left\{ \pmatrix f(x)& g(x) \\ 0 & 
f(x^2)\endpmatrix: f, g\in
k[x]\right\}\;\subset\; M_2(k[x]).
\end{equation}
Notice that this is the ring from \eqref{ex2.1}. The nilradical 
$N=N(R)$ is just the set of strictly upper triangular matrices. 
It is routine to check that $N$ is a free left $R/N$-module of rank one but a
free right $R/N$-module of rank two. Hence, by Theorem~\ref{thm1.1}, 
neither $R$ nor $S$ can have a
left noetherian  finite $\mathbb N$-filtrations.
\end{proof}

\begin{rem}\label{remark} 
The ring $R$ from \eqref{ex2.11}
clearly has Gelfand-Kirillov dimension one.
By inspecting the proof of Theorem~\ref{thm3.3}, this shows that 
there is no finite filtration of $R$ that induces a good filtration of
$N(R)$ as a left $R$-module.

Other examples of noetherian rings with no noetherian 
$\mathbb N$-filtration are given  in \cite{SZ}. However, those examples
require
that the ring in question have infinite Gelfand-Kirillov dimension 
and this is impossible for affine PI algebras (see \cite[Proposition
13.10.6]{MR}). 
\end{rem}

 It is readily checked that this ring $R$ does have a right noetherian,
finite filtration (see Corollary~\ref{cor4.fun}). By slightly modifying the
example one can get an example that ``works'' on both sides. 

\begin{exthm}\label{ex4.2} Let 
$$T=\left\{ \pmatrix f(x) & g(x) & h(x) \\
0 & f(x^2) & l(x) \\ 0&0& f(x)
\endpmatrix \; : \; f(x),\dots, l(x) \in k[x]\right\} + yM_3(k[x,y]).$$
Then $T$ is a affine  noetherian prime PI algebra such that, for every
$\mathbb N$-filtration, the associated graded ring is neither left nor
right noetherian. 
\end{exthm}

\begin{proof}
Notice that $T/(e_{13}T + e_{23}T) \cong R$, the ring 
from \eqref{ex2.11}, and so $T$ has no left noetherian finite 
filtration. Since $T\cong T^{op}$, the same is true 
on the right. 
\end{proof}

These examples can be easily modified into complete local rings and we 
give one example that is analogous to the ring $R$ from Example~\ref{ex4.1}.
Before giving the example, we need a definition. An ideal $I$ of a ring $R$ is
said to satisfy the {\it strong AR} property if the associated Rees ring 
$\Rees_IR=\bigoplus_{j\geq 0}I^j$ is noetherian. The significance of this
condition is that, if $I$ satisfies the strong AR property, then it also
satisfies the usual Artin-Rees (AR) property; indeed this is the standard way of
proving the latter condition in commutative algebra. The same idea has been
useful in noncommutative algebra (see, for example, \cite[\S2]{SW} 
and \cite{LV2}).
Given that the Jacobson radical of a semi-local 
noetherian PI algebra is automatically AR
(\cite[Theorems 3.1.13 and
 7.2.5]{J}), it is natural to ask if it is also strongly AR.

The next example shows that it does not. Notice that, although the ring in
question is just a completion of the ring $R$ from Example~\ref{ex4.1},
the Zariskian property fails on the opposite side.

\begin{exthm}\label{ex4.3.} 
Let  $$R=\left\{ \pmatrix f(x)& g(x) \\ 0 & f(x^2)\endpmatrix: f,g\in
k[[x]]\right\} + yM_2(k[[x,y]])\;\subset\; M_2(k[[x,y]]).$$ 
Then $R$ a prime noetherian complete local PI algebra over $k$.
Moreover,
$R$ has no weak adic filtration $\Gamma$ such that 
$\gr_\Gamma R$ is right noetherian. In particular, 
$J(R)$ does  not satisfy the strong AR property.
\end{exthm}

\begin{proof}
The proof that $R$ is prime noetherian PI algebra over $k$ is analogous to
that for Example~\ref{ex4.1} and is left to the reader. 
Clearly $R$ is a complete local ring. 

In order to  complete the proof it suffices, by Corollary~\ref{cor3.7},
  to work in a factor ring
and we chose: 
\begin{equation}\label{rone}
R'=R/yM_2(k[[x,y]]) \;=\;
\left\{\pmatrix f(x) & g(x) \\ 0 & f(x^2)\endpmatrix : f(x),g(x) \in
k[[x]]\right\}.
\end{equation}
    The nilradical $N$ of $R'$ is the set of strictly upper triangular matrices
    and  is free of rank $1$ as a left $R'/N$-module and free of
rank $2$ as a right $R'/N$-module.  
Now apply Corollary~\ref{cor3.7}. 
\end{proof}

Finally we justify a remark made after Theorem~\ref{thm3.3}: in that theorem
it is possible to have  a  filtration that is Zariskian on either the left or
the right.

\begin{cor}\label{cor4.fun}
Consider noetherian, PI 
algebras $R$ with ideals $I$ and $J_1\subseteq J_2$ such that
$J_2/J_1$ is a free left $R/I$-module of rank one and a free right $R/I$-module
of rank~$2$. 

Then, there exists an example $R_1$ of such a ring with a left Zariskian
filtration $\Gamma_1$ and 
an example $R_2$ of such a ring with a right Zariskian
filtration $\Gamma_2$. 
\end{cor}

\begin{proof}
The ring $R_2$ is just the ring $R$ from \eqref{ex2.1}, and so does satisfy the
hypotheses of the first paragraph. 
Set
\begin{equation}\label{alpha}
\alpha=\left(\matrix x&0\\0&x^2 \endmatrix\right)
\quad {\rm and }\quad \beta = \left(\matrix 0&1\\0&0 \endmatrix\right).
\end{equation}
It is easy to see that $R$ is generated by $\alpha$ and $\beta$ and hence
one has the standard $\mathbb N$-filtration $\Gamma_0=k$, 
$\Gamma_1=k+k\alpha + k\beta$ 
and $\Gamma_n = \Gamma_1^n$ for $n\geq 2$. 
Let $\aa $ and $\bb$ be the images of $\alpha$ and
$\beta$ in $\gr_{\Gamma} R$. Then, $\gr_{\Gamma} R$ is generated by 
$\aa$ and $\bb$ and 
satisfies the relations
 $\aa^2\bb=0=\bb\gr_\Gamma R\bb$.
Thus $\gr_{\Gamma} R$ is spanned by the elements 
$\{a^n, \bb\aa^n,\,
\aa\bb\aa^n\ : n\geq 0\}$. 
As such, $\gr_{\Gamma} R$ is a finitely generated right
$k[\aa]$-module, and so is right noetherian. 
As in the proof of Theorem~\ref{thm1.1}, this suffices to prove that 
$\Gamma$ is right Zariskian.

For $R_1$ we take the ring $R'$ from \eqref{rone}.
We define $\alpha$ and $\beta$ by \eqref{alpha} and observe that $\fm
=J(R') = R'\alpha + R'\beta$. Now use the $\fm$-adic
filtration $\Gamma$ defined by 
$\Gamma_i=R'$ if $i\geq 0$ but $\Gamma_i=\fm^{-i}$ if $i\leq 0$.
Then
$\beta\alpha=\alpha^2\beta\in \fm^3$. Hence in the associated graded ring
one finds that the images of these elements satisfy $\bb\aa=0=\bb^2$.
The argument of the last paragraph shows that 
$\gr_\Gamma R'$ is left noetherian and hence, by
\cite[Proposition~II.1.2.3]{LV2}, that  $\Gamma$ is left Zariskian.
\end{proof}

\section{A Dualizing Module }

As was remarked in the introduction, if an affine  noetherian PI algebra 
$R$ has a  finite 
noetherian $\mathbb N$-filtration, then $R$ has a dualizing 
complex. Thus, one can ask whether the ideas of the last section can be
used to
provide examples  of PI rings which do not have such a complex.
This appears not to be the case; in  this section  we check that
the ring $R$ from \eqref{ex2.1} does indeed have such a complex.

One advantage of this example is that we can work with modules rather than 
complexes, and so we define an $(R,R)$-bimodule $D$ to be a {\it dualizing
module} if (i) $D$ is finitely generated and of finite injective dimension on
both sides and (ii) the natural maps $R \rightarrow \End(D_R)$ 
and  $R^{\rm op} \rightarrow \End(_RD)$ are isomorphisms. 
A dualizing module viewed as a complex is a dualizing complex in the 
sense of Yekutieli \cite{Y}.

The way we find a dualizing module for the ring $R$ from \eqref{ex2.1}
is through the following observation:
Identify $C=k[x]$ with the
diagonal matrices in $R$ and set 
\begin{equation}\label{dualizing}
 D_1 = \Hom_{C}({}_{C}R, C)
\qquad{\rm and}\qquad
D_2 = \Hom_{C}({}R_{C}, C).\end{equation}
Thus, $D_1$ is an $(R,C)$-bimodule and 
$D_2$ is a $(C,R)$-bimodule.  The key to our construction is the following
easy lemma.

\begin{lem}\label{dual} 
Let $R\supseteq C$ be rings  such that $R$ is a finitely generated
projective $C$-module on both sides and $C$ is a commutative noetherian
algebra of finite injective dimension. 
Define modules $D_i$ by \eqref{dualizing}.
 Suppose that  
one has  ring isomorphisms
$\End_R(D_1)\cong R^{\rm op}$ and 
$\End_R(D_2)\cong R$ through which  $D_1\cong D_2$ as $R$-bimodules.
Then, $D=D_1$ is a dualizing module for $R$.
\end{lem}

\begin{proof} It is clear that $C$ is a dualizing module for itself.
 By \cite[Theorem~11.66]{Rot}
  one has natural isomorphisms: 
  \begin{equation}\label{rotman}
  \Ext^i_R({}_RN,{}_RD_1) =
\Ext_R^i({}_RN,\Hom_{C}({}_{C}R, C)) \cong \Ext^i_{C}({}_{C}N,C),
\end{equation}
 for any finitely generated left $R$-module $N$. This implies that
the injective dimension of $_RD_1$ is bounded by the injective
dimension of $_CC$. A similar assertion holds for $(D_2)_R$.
It follows from  \cite[Theorem 3.5]{SZ1} that $_RD_1$
and $(D_2)_R$ are finitely generated.  Therefore $D=D_1$ is a dualizing module.
\end{proof}

\begin{prop}\label{dual2}
 Let $R$ be the ring defined by \eqref{ex2.1} 
 and define modules $D_i$ by \eqref{dualizing}.  
Then, $D=D_1$ is a dualizing module for $R$.
\end{prop}

\begin{proof}
We check that the hypotheses of the lemma are satisfied for $C=k[x]$. 
Since $R$ is a free left  $k[x]$-module of rank two and a free right
$k[x]$-module of rank three, the proof is a routine computation which we will
only outline. One first checks that, under the natural module structures, 
$$D_1 \cong \frac{R\oplus R}{R(xe_{12},e_{12})}
\qquad{\rm and} \qquad
D_2 \cong R/e_{12}R,$$
where $e_{12}=\left(\smallmatrix 0&1\\0&0\endsmallmatrix\right)$ 
is the matrix unit. 
It follows that 
$$\End_R(D_2)\cong \II(e_{12}R)/e_{12}R\qquad{\rm where}\quad 
\II(aR)=\{\theta\in R : \theta aR\subseteq aR\}.$$
Computing this out one finds that 
\begin{align} \nonumber
\End_R(D_2) \cong& 
\left\{\pmatrix f(x^2)&g(x)\\0&f(x^4)\endpmatrix : f,g\in k[x]
\right\}\Big/e_{12}k[x^2]\\ \nonumber 
\cong&
\left\{\pmatrix f(x^2)&g(x^2)\\0&f(x^4)\endpmatrix : f,g\in k[x]
\right\} \cong R. \phantom{\pmatrix 1\\1\\1\endpmatrix}\nonumber
\end{align}

Finally, under this isomorphism $E=\End_R(D_2)\cong R$ 
one finds that the modules ${}_ED_2$ and ${}_RD_1$ become isomorphic,
 as is required to prove the
proposition.
\end{proof}

It would be interesting to know whether this proof can be extended to 
work for any ring $S$ that is finitely generated as a module over a 
commutative subring. 
The fact that the present proof depends upon the ``lucky'' isomorphism 
$\End_R(D_2)\cong R$ makes this seem unlikely. By using 
ideas from \cite{WZ}, it can at least be extended 
to Hopf algebras finitely generated as modules over commutative 
subalgebras. 

\bigskip
\centerline{\bf Acknowledgment}

\bigskip

The   authors thank Quanshui Wu for several conversations 
on the subject.

\ifx\undefined\bysame
\newcommand{\bysame}{\leavevmode\hbox to3em{\hrulefill}\,}
\fi

\end{document}